\documentclass[12pt]{article}
\usepackage{amsmath,enumerate, amsfonts, amssymb, amsthm}
\usepackage{eepic}
\usepackage{graphicx}
\def\UseSection{
        \numberwithin{equation}{section}
    \theoremstyle{plain}
        \newtheorem{theorem}    {Theorem}[section]
        \DefineTheorems 
}

\def\DefineTheorems{
    
    \newtheorem{lemma}      [theorem] {Lemma}
    
    \newtheorem{prop}       [theorem] {Proposition}
    
    \newtheorem{cor}        [theorem] {Corollary}

    \theoremstyle{definition}
    \newtheorem{defn}       [theorem] {Definition}

    \theoremstyle{definition}

}

\newcommand{\bt}   {\begin{theorem}}
\newcommand{\et}   {\end  {theorem}}
\newcommand{\bl}   {\begin{lemma}}
\newcommand{\el}   {\end  {lemma}}
\newcommand{\bp}   {\begin{prop}}
\newcommand{\ep}   {\end  {prop}}
\newcommand{\bc}   {\begin{cor}}
\newcommand{\ec}   {\end  {cor}}
\newcommand{\bd}   {\begin{defn}}
\newcommand{\ed}   {\end  {defn}}

\newcommand{\ba}   {\begin{array}}
\newcommand{\ea}   {\end  {array}}
\newcommand{\be}   {\begin{enumerate}}
\newcommand{\ee}   {\end  {enumerate}}
\newcommand{\bi}   {\begin{itemize}}
\newcommand{\ei}   {\end  {itemize}}

\def\eq#1\en{\begin{equation}#1\end{equation}}
\def\eqsplit#1\ensplit{
    \begin{equation}\begin{split}#1\end{split}\end{equation}
    }
\def\eqalign#1\enalign{
    \begin{align}#1\end{align}
    }
\def\eqmul#1\enmul{
    \begin{multline}#1\end{multline}
    }
\newcommand{\eqarrstar} {\begin{eqnarray*}}
\newcommand{\enarrstar} {\end{eqnarray*}}
\newcommand{\eqarray}   {\begin{eqnarray}}
\newcommand{\enarray}   {\end{eqnarray}}

\newcommand{\lbeq}[1]  {\label{e:#1}}
\newcommand{\refeq}[1] {\eqref{e:#1}}    

\makeatletter
\newcommand{\labelcounter}[2]{{%
    \stepcounter{#1}
    \protected@write\@auxout{}%
    {\string\newlabel{#2}{{\csname the#1\endcsname}{\thepage}}}%
    {\ref{#2}}
    }}
\makeatother
\newcommand{\sss}   { \scriptscriptstyle }

\newcommand{\Ebold} {{\mathbb E}}

\newcommand{\Nbold} {{\mathbb N}}

\newcommand{\Pbold} {{\mathbb P}}

\newcommand{\Rbold} {{\mathbb R}}

\newcommand{\Dcal}   {\mathcal{D}}

\newcommand{\spose}[1] {{\hbox to 0pt{#1\hss}} }
\newcommand{\ltapprox} {\mathrel{\spose{\lower 3pt\hbox{$\mathchar"218$}}
 \raise 2.0pt\hbox{$\mathchar"13C$}}}
\newcommand{\gtapprox} {\mathrel{\spose{\lower 3pt\hbox{$\mathchar"218$}}
 \raise 2.0pt\hbox{$\mathchar"13E$}}}

\UseSection  
\newcommand{\E}{\mathbb{E}}

\newcommand{\Var}{\operatorname{Var}}

\newcommand{\eps}{\varepsilon}
\newcommand{\e}{\operatorname e}

\newcommand{\nn}{\nonumber}

\def\1{{\mathchoice {1\mskip-4mu\mathrm l} 
{1\mskip-4mu\mathrm l}
{1\mskip-4.5mu\mathrm l} {1\mskip-5mu\mathrm l}}}

\newcommand{\ITk}{I^{\sss ({\rm T})}_{{\rm out}_k}}
\newcommand{\IBk}{I^{\sss ({\rm B})}_{{\rm out}_k}}
\newcommand{\su}{\sigma_u}
\newcommand{\ub}{\bar{I}_u}
\newcommand{\INLm}{{{\rm INL}}_{{\rm max}}}
\newcommand{\eqd}{\stackrel{\Dcal}{=}} 
\newcommand{\Ilsb}{{I}_{{\rm lsb}}}

\newcommand{\Nfat}{\mathbf{N}}
\newcommand{\Zfat}{\mathbf{Z}}
\newcommand{\n}{\mathbf{n}}
\newcommand{\M}{\overline M}

\def\1{{\mathchoice {1\mskip-4mu\mathrm l} 
{1\mskip-4mu\mathrm l}
{1\mskip-4.5mu\mathrm l} {1\mskip-5mu\mathrm l}}}

\title{Functionals of Brownian bridges \\arising in the current mismatch in D/A-converters}

\author{Markus Heydenreich$^1$ \and
        Remco van der Hofstad$^1$
        \and
        Georgi Radulov$^2$
    }

\begin{document}
\maketitle
\footnotetext[1]{Department of Mathematics and
        Computer Science, Eindhoven University of Technology,
        5600~MB Eindhoven, The~Netherlands.
        {\tt m.o.heydenreich@tue.nl, r.w.v.d.hofstad@tue.nl}}
\footnotetext[2]{Department of Electrical Engineering,
        Eindhoven University of Technology, EH 5.15,
        5600~MB Eindhoven, The~Netherlands.
        {\tt g.radulov@tue.nl}}

\begin{abstract}
\noindent
Digital-to-analog converters (DAC) transform signals from the abstract
digital domain to the real analog world. In many applications, DAC's
play a crucial role.

Due to variability in the production, various errors
arise that influence the performance of the DAC. We focus on the
current errors, which describe the fluctuations in the currents
of the various unit current elements in the DAC. A key performance
measure of the DAC is the {\it Integrated Non-linearity} (INL),
which we study in this paper.

There are several DAC architectures. The most widely used architectures are
the thermometer, the binary and the segmented architectures.
We study the two extreme architectures, namely, the thermometer and the binary
architectures. We assume that the current errors are i.i.d.\ normally distributed,
and reformulate the INL as a functional of a Brownian bridge. We then
proceed by investigating these functionals. For the thermometer case, the functional
is the maximal absolute value of the Brownian bridge, which has been
investigated in the literature. For the binary case, we investigate
properties of the functional, such as its mean, variance and density.
\end{abstract}

\section{Current mismatch in digital-to-analog converters}

Digital-to-Analog converters (DAC) transform signals from the abstract digital domain to the real analog world.
For many applications, this conversion enables the usage of the computational power of robust digital electronics.
For example, digital audio and video, digital control, and telecommunications are fields that require digital-to-analog conversion.
The advantageous intelligence of the applications in these fields is implemented with digital logic, e.g. microprocessors, and used via DA conversion in the real analog world.
However, the DAC errors, at the end of the application chain, may decrease the performance of the whole system.
Therefore, predicting and controlling these errors is crucial.
This requirement is further emphasized in the highly integrated mixed-signal Systems-on-a-Chip (SoC), as we will explain now.

Nowadays, the SoC system solutions often integrate DAC functionality together with the digital logic for cost effectiveness.
This requires optimal usage of the DAC resources while keeping the errors of the DAC within specified margins, primarily because of the two following reasons.
Firstly, the price of the mixed-signal SoC includes the price of the DAC resources even when customers are not interested in the DAC functionality.
Secondly, if the DAC does not comply with its specifications, then the entire SoC chip should be discarded.
By a careful design, the DAC errors mainly arise from the uncertainty in the manufacturing process and, hence, they are random.
Therefore, statistical rules are used to predict the overall performance for high-volume chip production.
Knowledge is required that accurately links the DAC resources with the DAC error margins.
An important example of such a relationship is the dependence of the Digital-to-Analog (D/A) conversion accuracy on the DAC area, i.e.\ the silicon area of the chip that is used for the DAC.
Higher conversion accuracy is achieved for larger chip areas \cite{Bastos01, Pelgrom89}.
However, too large areas introduce additional problems which may degrade the conversion accuracy.
Thus, precise knowledge and understanding of the relationship between accuracy and chip area is crucial, particularly for high-volume chip systems including DAC.
For a general introduction to DAC converters we refer to \cite{Bastos01, Jesper01,Plassche94}, whereas Razavi \cite{Razavi01} also focusses technical aspects.

The D/A conversion is carried out by switching certain analog quantities, such as voltages, currents or charges, ON or OFF.
For the sake of simplicity, and without loss of generality, this paper will assume current as a basic analog quantity. The switching process is controlled by the digital input signal $w\in\{0,1\}^N$, where $N$ is the length of the binary input signal.
If all combinations of 0's and 1's between the  digital bits produce valid input signals, i.e., $N=\log_2(n+1)$,
then the coding is called \emph{binary} and $N$ is called the \emph{resolution} of the DAC.
A DAC that uses binary coding to control its current quantities is called binary DAC.
The switched ON current quantities $I_{u_i}$ are summed to construct the analog
output signal current $I_\text{out}$. We will assume that the current quantities are
{\it random} and that $\{I_{u_j}\}_{j=1}^n$ are i.i.d.\ random variables.
The sum of all current quantities
    \eq
    I^{}_{\text{out}_\text{max}}=\sum_{j=1}^n I_{u_j}
    \en
is associated to the maximal digital input word $w_n=(1,1,\dots,1)$, and is called
the full-scale current of the DAC.
The smallest meaningful difference in the analog output, defined as the output
{\it least-significant bit} (LSB), is defined as the full-scale current
divided by the number of digital input codes, i.e.,\
    \eq\lbeq{defIlsb}
    I_{\text{lsb}}:=\frac{I^{}_{\text{out}_\text{max}}}{n}.
    \en
For general DAC, errors can be classified as static or dynamic.
We will focus on the \emph{static} errors, which we now introduce.

For every digital input word $w_k\in\{0,1\}^N$, there is an analog output value
$I^{}_{\text{out}_k}$. The code words $w_k\in\{0,1\}^N$ will be assumed to be ordered.
The difference between two adjacent output values would ideally be $\ub=\E[I_{u_j}]$, but in practice it deviates due to mismatch errors coming from the uncertainty of the manufacturing process, i.e.\
    \eq
    I^{}_{\text{out}_k}-I^{}_{\text{out}_{k-1}}=\ub+\left(\Ilsb-\ub\right) +\text{DNL}_k\cdot\Ilsb,
    \qquad k=1,\dots,n.
    \en
Here, $\Ilsb-\ub$ represents the {\it linear error} which is independent of $k$, and
    \eq\label{eqDefDNLk}
    \text{DNL}_k=\frac{I_{\text{out}_k}-I_{\text{out}_{k-1}}-I_{\text{lsb}}}{I_{\text{lsb}}}.
    \en
is the \emph{differential non-linearity} in $\Ilsb$-scale.
As the DAC are required to be \emph{linear} devices, the non-linear errors are the main concern.
The {\it integrated non-linearity} (INL) measures the non-linearity
of the whole DAC transfer characteristic, i.e., the cumulated individual non-linear errors.
This paper concentrates on the INL, as defined for a code $w_k$ by
    \eq
    \label{eqDefINLk}
    \text{INL}_k:=\sum_{i=1}^k\text{DNL}_i=\frac{I_{\text{out}_k}-k\cdot I_{\text{lsb}}}{I_{\text{lsb}}},
    \qquad k=1,\dots,n,
    \en
This definition excludes the linear errors by forcing $\text{INL}_0=\text{INL}_{n}=0$ and normalizes $\text{INL}_k$ to LSB-scale.
The maximal absolute INL over all codes is given by
    \eq
    {\text{INL}}_{\text{max}}:=\max_{k=0,\dots,n}\left|\text{INL}_k\right|.
    \en
The statistic $\INLm$ is an important specification, because it indicates how linear the DAC is.
Another important figure is $\text{Yield}_\text{INL}$, which indicates how many of the manufactured chips have ${\text{INL}}_{\text{max}}$ under certain specified threshold.
This is to say, if a DAC manufacturer guarantees certain linearity, then
$\text{Yield}_\text{INL}$ describes what proportion of the produced chips fall
within these specifications.

INL is the most popular DAC specification, see e.g.\ \cite{Jesper01,Plassche94}.
Its practical importance is very high in all DAC application fields.
This general definition of INL is valid for all DAC architectures, and all
DAC resolutions. The most commonly used architectures are the binary, segmented,
and thermometer architectures. The segmented architecture interpolates between
the binary and the thermometer architectures.
In this paper, we consider the two extreme cases of \emph{thermometer} and
\emph{binary} DAC architectures. We now describe these DAC architectures.

For a \emph{thermometer} DAC, $I_{\text{out}_k}=\ITk$ is given by
\eq\lbeq{defIoutThermo}
    \ITk=\sum_{j=1}^{k} I_{u_{j}},\qquad k=0,\dots,n.
\en
\begin{figure}
    \centering
    \includegraphics[width=10cm]{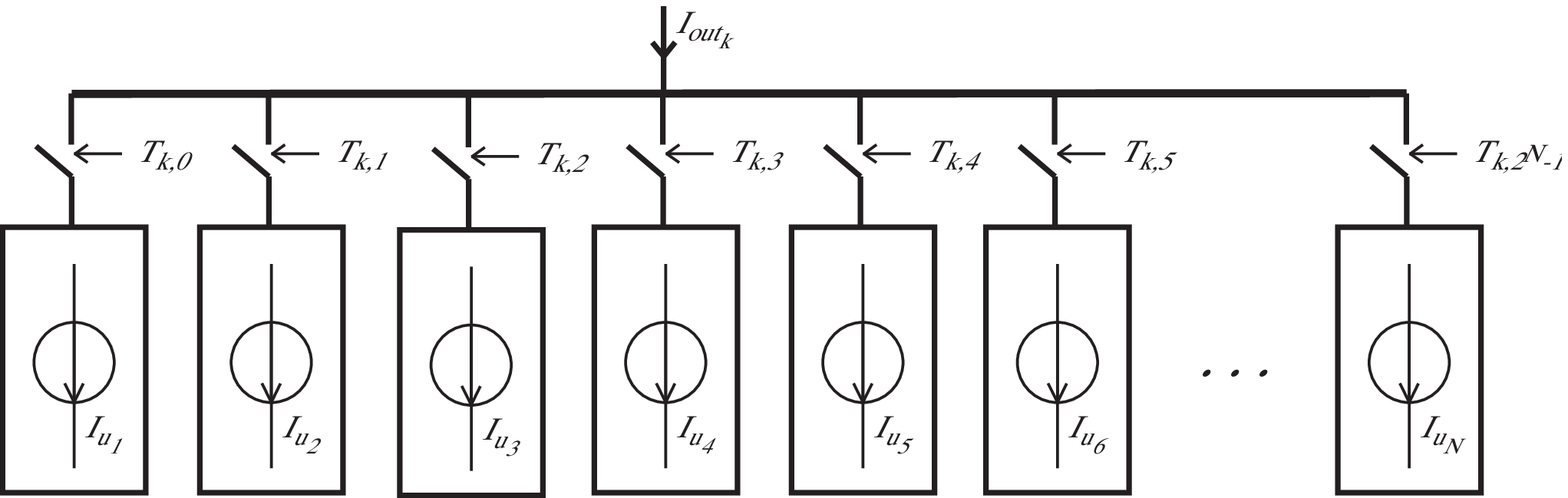}
    \caption{Thermometer DAC. The code word $w_k$
    corresponds to switching $T_{k,i}$ ON for $i\leq k$ and
    switching $T_{k,i}$ OFF for $i>k$.}
    \includegraphics[width=10cm]{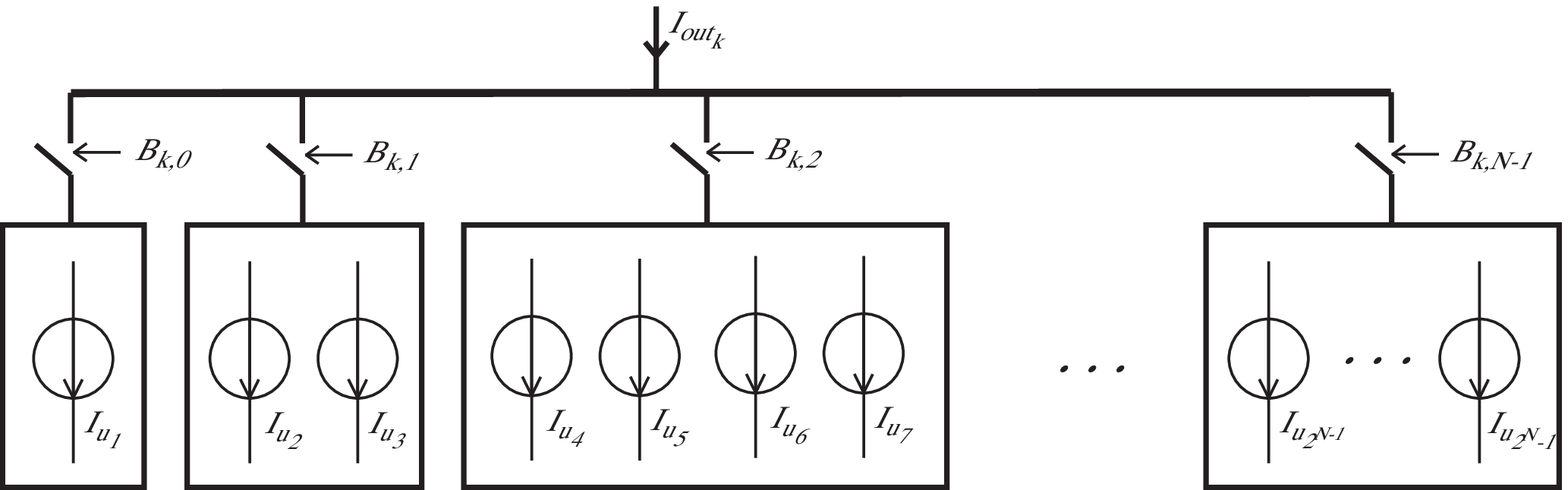}
    \caption{Binary DAC. The code word $w_k$
    corresponds to switching those $i$ for which $B_{k,i}=1$ ON, whereas
    the $i$ for which $B_{k,i}=0$ are switched OFF.
    The matrix $B$ is given in \refeq{def-B}.}
    \includegraphics[width=10cm]{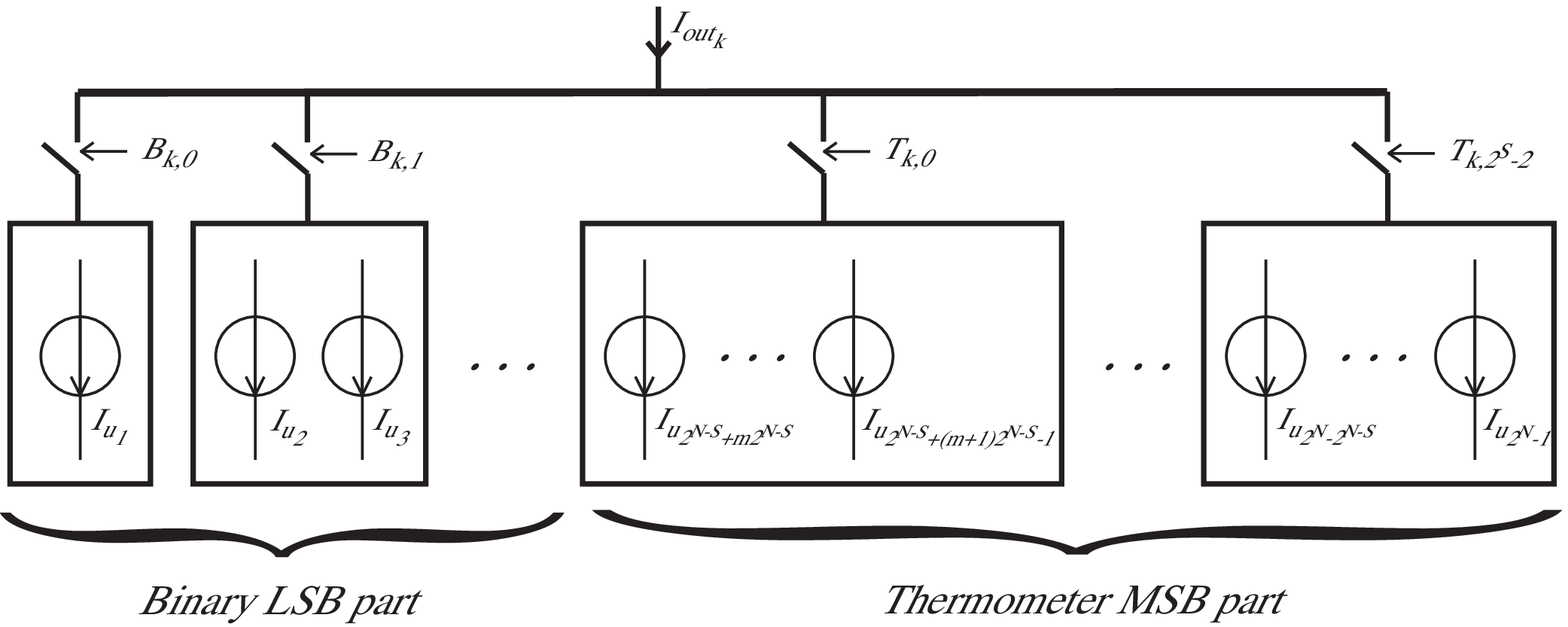}
    \caption{Segmented DAC.}
    \label{figSegmented}
\end{figure}

For a \emph{binary} DAC, on the other hand, $I_{\text{out}_k}=\IBk$ is given by
\eq\lbeq{defIoutBin}
    \IBk=\sum_{m=1}^{N} B_{k,m} \sum_{j=2^{m-1}}^{2^{m}-1}I_{u_{j}},\qquad k=0,\dots,n,
\en
where the {\it switching matrix} $B$ is an ${n\times N}$-matrix given by
\eq
\lbeq{def-B}
    B=
  \left(%
  \begin{array}{cccccc}
  0 & 0 & 0 & 0 & \cdots & 0 \\
  1 & 0 & 0 & 0 & \cdots & 0 \\
  0 & 1 & 0 & 0 & \cdots & 0 \\
  1 & 1 & 0 & 0 & \cdots & 0 \\
  0 & 0 & 1 & 0 & \cdots & 0 \\
  \vdots & \vdots & \vdots & \vdots & \ddots & \vdots \\
  1 & 1 & 1 & 1 & \cdots & 1 \\
  \end{array}%
  \right).
\en
The matrix $B$ is constructed by writing the first $n$ integers in
(reversed order) binary coding into the rows of $B$.
Note that the columns of the matrix represent the DAC bits, i.e., the
switches of the grouped current sources.
The 0's represent switched OFF currents, and the 1's represent
switched ON currents.
The most left column gives the switches for the least-significant bit current
$I^{\sss (B)}_{\text{out}_1}=I_{u_1}$, while the most right
column gives the switches for the {\it most-significant bit} (MSB) current
$I^{\sss (B)}_{\text{out}_N}=\sum_{j=2^{N-1}}^{2^N-1}I_{u_j}$.
Note that
    \eq
    \lbeq{Ioutmax}
    I^{\sss (B)}_{\text{out}_n}=\sum_{m=1}^{N} B_{n,m} \sum_{j=2^{m-1}}^{2^{m}-1}I_{u_{j}}=\sum_{j=1}^n I_{u_j}
    =n\Ilsb,
    \en
so that $\text{INL}_n=0.$

In practice, the most popular DAC architecture is the segmented one. The segmented DAC architectures implement one part of the input digital bits in a binary way and the other part in a thermometer way. That is how the advantages of both the binary and the thermometer architectures can be combined.

Figure \ref{figSegmented} shows a segmented DAC architecture. The LSB part is implemented in a binary way and the MSB part in a thermometer way.
The output current for a given level of segmentation $S$ is expressed by
    \eq
    \lbeq{defIoutSegm}
    {I^{\sss (\text{S})}_{\text{out}_k}}
    =\sum_{m=1}^{N-S}B_{k,m}\sum_{j=2^{m-1}}^{2^{m}-1}I_{u_j}+\sum_{m=1}^{\lfloor k/2^{N-S}\rfloor}\sum_{j=m2^{N-S}}^{(m+1)2^{N-S}-1}I_{u_j}.
    \en
The parameter $S\in\{0,\dots,N\}$ determines the interpolation between
the binary and the thermometer parts of the DAC architecture.
Indeed, for $S=0$, the segmented DAC architecture is transformed to a fully binary architecture.
On the other hand, for $S=N$, the segmented DAC architecture is transformed to a fully thermometer architecture.
Note that in both extreme cases the same unit current sources $I_{u_j}$ are used.
Thus, the performance difference is only due to the way the $I_{u_j}$ are combined to construct the output current $I_{\text{out}_k}$.
In the binary DAC, the unit currents are first grouped and then switched ON or OFF, while in the thermometer DAC, the unit currents are individually switched ON or OFF.
More detailed discussion on DAC architectures can be found in the literature \cite{Jesper01, Plassche94, Wikner01}.

As discussed, due to manufacturing related mismatch, the currents $I_{u_{j}}$ always deviate from their designed values, i.e., $I_{u_{j}}=\ub+\eps_j$.
The mean value $\ub$ is chosen by the DAC designer, whereas $\eps_j$ is the random error due to the manufacturing process which we model as an i.i.d.\ sequence of normal random variable with zero mean and variance $\su^2$.
Nevertheless, our results remain valid when the $\{\eps_j\}$ are i.i.d.\ distributed with sufficiently many moments.
The ratio $\su/\ub$ is known as the {\it relative current matching}, i.e., the unit current matching.
The relative matching determines the required transistor area for the particular manufacturing process.
The smaller $\su/\ub$, the more accurate $I_{u_j}$, but the larger the required area.
For more details on the dependence between relative matching and transistor area we refer to \cite{Pelgrom89}.

Once $\su/\ub$ is specified, the required transistor area can be calculated.
However, the relationship between $\su/\ub$ and DAC non-linearity, which is crucial to determine the proportion of chips complying to specifications, has never been determined analytically.
So far, DAC engineers have used either Monte Carlo simulations or approximations.

The Monte Carlo simulations of a DAC model produce empirical results to suggest some design specifications $\su/\ub$. Although not accurate, this approach is very practical, see \cite{Conroy89}.
A problem arises for the design of high-resolution DAC, for which $N$ is large.
The complexity of the DAC model increases by a factor of two for every additional bit in resolution. For example, for $N=14$, $n=16383$ unit elements have to be simulated. Therefore, Monte Carlo simulations are not practical for higher resolutions, because they become complex and slow.

On the other hand, a number of analytical approximations can be found in the literature. The analytical attempts to describe INL, and in particular $\text{INL}_\text{max}$, started in 1986. The approach in \cite{Laksch86} disregards the correlation between the DAC outputs for different input codes.
Bastos \cite{Bastos01} proposes a much simpler formula, which considers only the deviation of the transfer characteristic at the mid-scale DAC output, which can be a rough, though simple, estimation of the $\text{INL}_\text{max}$.
Another approximation was given by van den Bosch et al.\ \cite{Bosch04} by assuming that, if the DAC static transfer characteristic at any code INL has error equal to the target value, e.g., $\text{INL}_k=0.5\,\Ilsb$, then there should be a 50\% chance that ultimately $\text{INL}_\text{max}$ is smaller than the target value, i.e.,  $\text{INL}_k\le0.5\,\Ilsb$ for all $k$. The major approximation inaccuracy is in the probability that both the positive and negative INL limits are reached for the same DAC sample, e.g., $\text{INL}_k<-0.5\,\Ilsb$ is disregarded. Though this approach derives a convenient normal distribution for $\text{INL}_\text{max}$, it is inaccurate for higher resolutions, as we show in more detail in this
paper. In general, approximations lead to transistor over-design, i.e., a transistor area that is too large.
In this paper, on the other hand, we shall present an exact analytical formula, for which no approximation
is necessary.

Due to the lack of exact analytical formulation of INL and the high complexity of DAC model simulations, the statistics used in industry for a high volume chip production are hard to predict.
Here we think of the statistics $\text{Yield}_\text{INL}$, the $\text{INL}_\text{max}$ distribution, $\text{INL}_\text{max}$ deviation and mean.
Furthermore, the advantages of some redundancy based approaches relying on the DAC statistical INL properties cannot be theoretically estimated, see e.g.\ \cite{Radulov06}. Finally, up to now the main DAC architectures, i.e.\ binary, thermometer, and segmented, cannot be distinguished with respect to their static linearity properties, so they are wrongly considered identical \cite{Bastos01,Bosch04}. One conclusion from our results is that the INL
for binary and thermometer architectures are {\it different}.

Implications of the results derived in this paper for the field of
DAC's can be found in a second paper \cite{RaduHeydHofsHegtRoer06}.
A comparison with the results in \cite{Bosch04,Conroy89,Laksch86} is summarized in \cite[Table 1]{RaduHeydHofsHegtRoer06}.

\section{Thermometer coding: Maximum of a Brownian bridge}
For the thermometer coding, we can describe the INL as functional of a Brownian bridge
as follows.
\bt[$\INLm$ for the \emph{thermometer coding}]\label{thm-ThermoResults}
    As $n\to\infty$,
    \eq
    \lbeq{INLconv}
        \frac{\Ilsb}{\su\sqrt{n}}\INLm\longrightarrow X,
    \en
    in distribution, in $L^1$ and $L^2$, where the limit $X$ is characterized by
    \eq
        X=\max_{t\in[0,1]}\left|B_{t}\right|
    \en
    for a Brownian bridge $\{B_s\}_{s\in[0,1]}$, and
    \eq\lbeq{MomentsOfX}
        \E[X]=\frac{\sqrt{2\pi}}{2}\ln2\approx0.86875,
        \qquad
        \Var(X)=\frac{\pi^2}{12}-\frac{\pi}{2}(\ln2)^2\approx0.0677732,
        \qquad
    \en
    and
    \eq\lbeq{UpperTailX}
        \Pbold(X\le x)=1+2\sum_{k=1}^\infty(-1)^k\e^{-2k^2x^2},    \qquad x>0.
    \en

\et

Note that in \refeq{INLconv}, we multiply by $\Ilsb$ rather than by
$\ub$. For the convergence in distribution, this makes no difference
what so ever, since $\Ilsb$ converges to $\ub$ a.s.\ by the strong
law of large numbers. However, the convergence in $L^1$ and $L^2$
{\it fails} if we multiply by $\ub$, since in this case,
the expected value of ${\rm INL}_k=({I_{\text{out}_k}-k\cdot I_{\text{lsb}}})/{I_{\text{lsb}}}$
is not defined, whereas $|{\rm INL}_k|$ has infinite mean.

We recall that a Brownian bridge $\{B_s\}_{s\in[0,1]}$ is a Markov
process on $[0,1]$ that is obtained from a Wiener process or Brownian motion
in either of the following two ways:
\begin{enumerate}
\item[(B1)]
    \begin{equation*}
    B_s=W_s-sW_1,
    \end{equation*}
where $\{W_s\}_{s\in[0,1]}$ is a Wiener process;
\item[(B2)]\label{B2}
    \begin{equation*}
    B_s=W_s',
    \end{equation*}
where $\{W_s'\}_{s\in[0,1]}$ is a Wiener process conditioned on $W_1'=0$.
\end{enumerate}
For an introduction to Wiener processes and Brownian bridges we refer to \cite{Grimmett92}.
For the equivalence of (B1) and (B2) cf.\ \cite[pp.\ 83--85]{Billingsley68}.

\bp\label{thermo}
For the thermometer coding,
    \eq\lbeq{thermo1}
    \max_{k=0,\dots,n}\left|\ITk-k\, \Ilsb\right|
    \eqd \su\sqrt{n}\max_{k=1,\dots,n}\left|B_{k/n}\right|,
    \en
where $\{B_t\}_{t\in[0,1]}$ is a Brownian Bridge process and $\eqd$ denotes equality in distribution.
\ep
\proof
Since $\{I_{u_j}\}_{j=1,\dots,n}$ is a family of i.i.d.\ normally
distributed random variables with mean $\ub$ and variance $\su^2$,
we have that
    \eq
    \ITk-k\, \Ilsb=\sum_{j=1}^k\left(I_{u_j}-k\,\ub\right)
    \en
is normally distributed with mean $0$ and variance $k\su^2$.
Hence for a Brownian motion $\{W_t\}_{t\ge0}$,
    \eq\lbeq{thermoscaling1}
    \ITk-k\,\Ilsb\eqd\su W_k\eqd\su\sqrt{n}\, W_{k/n},
    \en
where we used Brownian scaling in the last distributional equality.
By substituting $k=n$ and using \refeq{defIlsb} we obtain $n(\Ilsb-n\ub)\eqd\su\sqrt{n} W_1$.
Combined with \refeq{thermoscaling1}, this yields
    \eq\lbeq{thermoscaling2}
    \ITk-k\,\Ilsb\eqd\su\sqrt{n} \left(W_{k/n}-\frac kn W_1\right)\eqd\su\sqrt{n}\, B_{k/n}
    \en
for a Brownian bridge process $\{B_s\}_{s\in[0,1]}$, where we used (B1) above.
After taking absolute values on both sides of \refeq{thermoscaling2} and the maximum over $k=1,\dots,n$, we obtain \refeq{thermo1}.
\qed
\vspace{0.5cm}

The maximum over the discrete time points $k/n$, $k=1,\dots,n$, in \refeq{thermo1}
can be replaced by the supremum over the whole interval $[0,1]$ by using the following lemma:
\bl\label{ConvergenceMaxB}
    For $C>4$,
    \eq\lbeq{thermo7}
    \Pbold\left(\Big|\max_{k=1,\dots,n}\left|B_{k/n}\right|-\max_{t\in[0,1]}\left|B_{t}\right|\Big|
    \ge C\sqrt{\frac{\log n}{n}}\right)
    \le4n^{1-C^2/8}+2n^{-C^2n/8}.
    \en
    In particular, $\max_{k=1,\dots,n}\left|B_{k/n}\right|$ converges to $\max_{t\in[0,1]}\left|B_{t}\right|$ in distribution as $n\to\infty$.
    Moreover, $L^1$- and $L^2$-convergence holds.
\el

The proof of Lemma \ref{ConvergenceMaxB} is deferred to the appendix.
We now use Lemma \ref{ConvergenceMaxB} to complete the proof of
Theorem \ref{thm-ThermoResults}:

\proof[Proof of Theorem \ref{thm-ThermoResults}.]
The convergence in \refeq{INLconv} follows from Lemmas \ref{thermo} and \ref{ConvergenceMaxB}.
For the probability of the upper tail of $X$ \refeq{UpperTailX} we refer to
\cite[(11.39)]{Billingsley68}, the computation of mean and variance of $X$
is straightforward integration. For example, we have that
    \eq
    \Ebold[X]=\int_0^{\infty} \Pbold(X>x)dx
    =2\int_0^{\infty} \sum_{k=1}^\infty(-1)^{k-1}\e^{-2k^2x^2}dx
    =2\sum_{k=1}^\infty(-1)^{k-1} \frac{\sqrt{2\pi}}{4k}
    =\frac{\sqrt{2\pi}}{2}\ln{2}.
    \en
The interchange of summation and integral can be justified by
looking at $X_{\varepsilon}=X\vee \varepsilon$.
\qed
\vspace{.5cm}

\section{Binary coding: Block increments of a Brownian bridge}
\subsection{Results and overview proof}
In this section we prove the following theorem characterizing the binary coding statistic. \bt[$\INLm$ for the \emph{binary coding}]\label{thm-BinaryResults}
    As $n\to\infty$,
    \eq\lbeq{limitM}
        \frac{\Ilsb}{\su\sqrt{n}}\INLm\longrightarrow M,
    \en
    in distribution, in $L^1$ and $L^2$, where the limit $M$ is characterized by
    \eq\lbeq{formulaM}
        M=\frac12\sum_{l=1}^\infty2^{-(l+1)/2}\left|\sum_{j=1}^{l-1}2^{-(l-j)}Z_j-Z_l\right|
    \en
    for a family $\{Z_l\}_{l=1}^\infty$ of i.i.d.\ standard normal random variables.
    The expectation of $M$ is given by
    \eq\lbeq{expecM}
        \E[M]=(2\pi)^{-1/2}\sum_{l=1}^\infty\big(2^{-l}-2^{-2l}\big)^{1/2}
        \approx0.839792,
    \en
    and the variance $\Var(M)$ is computed explicitly in \refeq{formulaVar} below and can be approximated as
    \eq\lbeq{varianceOfM}
        \Var(M)\approx0.08007 .
    \en
\et

Note that in \refeq{limitM}, we multiply by $\Ilsb$ rather than by
$\ub$. As explained for Theorem \ref{thm-ThermoResults}, this makes no difference for the convergence in distribution, though the convergence in $L^1$ and $L^2$ fails.

The proof of Theorem \ref{thm-BinaryResults} is organized as follows.
In Lemma \ref{MNconvM} we prove the convergence in \refeq{limitM},
where the limit is characterized in terms of increments of Brownian bridges.
After this, Proposition \ref{thm-1} shows that the weak limit $M$ can be
expressed as the weighted sum of standard normal variables, which proves \refeq{formulaM}.
Lemma \ref{MomentsOfM} states the expression for mean and variance of $M$.
Finally, we give an approximation to the density of $M$ in Section \ref{sec-densityM}.

\subsection{A Brownian bridge representation of the binary $\INLm$}
Our aim is to derive an expression for $\INLm$ for the binary coding.
First we express the non-linearity of the current steering DAC in terms of a functional of a Brownian bridge.
\bl\label{MNformula}
    \eq\label{eqM}
    \max_{k=1,\dots,n}\left|\IBk-k\Ilsb\right|
    \eqd \frac \su2\sqrt{n} \sum_{m=1}^{N}\left|B_{\frac{2^m-1}n}-B_\frac{2^{m-1}-1}n\right|,
    \en
where 
$\{B_s\}_{s=0}^1$ is a Brownian bridge.
\el
\proof
Let $\{W_t\}_{t\in[0,\infty)}$ be a Wiener process, then
\eq
    I_{u_j}-\ub\eqd\su\sqrt{n}\left(W_{j/n}-W_{(j-1)/n}\right),\qquad j=1,\dots,n,
\en
where $\eqd$ represents equality in distribution. We will further write $=$ rather than
$\eqd$, because we are interested in the distribution only.
Furthermore,
\eq
    I_{\text{lsb}}-\ub
    =\frac1n\sum_{j=1}^n\left(I_{u_j}-\ub\right)
    =\frac1n\,\su W_n
    =\frac{1}{\sqrt{n}}\,\su W_1,
\en
where we recall that the $I_{u_j}$ are i.i.d.\ $\mathcal{N}(\ub,\su)$-distributed.

We want to calculate $I^{(B)}_{\text{out}_k}-kI_\text{lsb}$ for $k$ being a power of 2 first. The advantage is that, for $k=2^{m-1}$, only the $m$th block
$\{I_{u_j}|j={2^{m-1}},\dots,{2^{m}-1}\}$
contributes:
\begin{eqnarray*}
  I^{(B)}_{\text{out}_{2^{m-1}}}-2^{m-1}I_\text{lsb}
  &=&\left(I^{(B)}_{\text{out}_{2^{m-1}}}-2^{m-1}\ub\right)+2^{m-1}\left(\ub-I_\text{lsb}\right)\\
  &=&\sum_{j=2^{m-1}}^{2^{m}-1}\left(I_{u_j}-\ub\right)+2^{m-1}\left(\ub-I_\text{lsb}\right)\\
  &=&\su\sqrt{n}\left(\Big(W_{\frac{2^{m}-1}n}-\frac{2^{m}-1}nW_1\Big)-\Big(W_{\frac{2^{m-1}-1}n}-\frac{2^{m-1}-1}{n}W_1\Big)\right)\\
  &=& \su\sqrt{n}\left(B_{\frac{2^{m}-1}n}-B_{\frac{2^{m-1}-1}n}\right),
\end{eqnarray*}
where, in the last line, we used the representation (B1) of Brownian bridges.

For calculating $\max_{k=1,\dots,n}\big|I^{(B)}_{\text{out}_{k}}-k\,I_\text{lsb}\big|$, we need to take the maximum over every configuration of contributing blocks, hence
\begin{eqnarray}\label{eqM1}
\max_{k=1,\dots,n}\left|I^{(B)}_{\text{out}_{k}}-k\,I_\text{lsb}\right|
&=&\su\sqrt{n}\,\underbrace{\max_{\mathfrak{I}\subset\left\{1,...,N\right\}}\left|\sum_{m\in\mathfrak{I}}B_{\frac{2^m-1}n}-B_\frac{2^{m-1}-1}n\right|}_{:=M_N}.
\end{eqnarray}
We denote by $\mathfrak{I}^*$ the subset of $\{1,\dots,N\}$ for which the maximum in (\ref{eqM1}) is achieved,
and use the abbreviation
    \eq
    \hat B_{n,m}:=B_{\frac{2^m-1}n}-B_{\frac{2^{m-1}-1}n},\qquad m=1,\dots,N,
    \en
for the increment of the Brownian bridge on the interval $\big[\frac{2^{m-1}-1}n,\frac{2^m-1}n\big]$.
Without loss of generality, we may assume that
$\sum_{m\in\mathfrak{I}^*}\hat B_{n,m}$
is positive, otherwise the same argument for $-B$ holds. Clearly,
\eq\label{eqM2}
    m\in\mathfrak{I}^*\qquad\Leftrightarrow\qquad \hat B_{n,m}\ge0.
\en
Furthermore,
\eq\label{eqM3}
0=B_1-B_0
=\sum_{m=1}^N\hat B_{n,m}
=\sum_{m=1}^N\hat B_{n,m}\,\1_{\hat B_{n,m}\ge0}+\sum_{m=1}^N\hat B_{n,m}\,\1_{\hat B_{n,m}\le0}.
\en

Using (\ref{eqM2}) in the first line and (\ref{eqM3}) in the second, we obtain
\begin{eqnarray*}
  M_N
  &=& \sum_{m=1}^N\hat B_{n,m}\,\1_{\hat B_{n,m}\ge0} \\
  &=& \frac12\sum_{m=1}^N\hat B_{n,m}\,\1_{\hat B_{n,m}\ge0} +\frac12\sum_{m=1}^N(-\hat B_{n,m})\,\1_{\hat B_{n,m}\le0} \\
  &=& \frac12\sum_{m=1}^N\left|\hat B_{n,m}\right|,
\end{eqnarray*}
which is (\ref{eqM}).
\qed\vspace{.5cm}

We write
    \eq\lbeq{defMN}
    M_N=\frac1{\su\sqrt{n}}\max_{k=1,\dots,n}\left|I^{(B)}_{\text{out}_{k}}-k\,I_\text{lsb}\right|
    \en
as in (\ref{eqM1}), and define
    \eq\lbeq{defM}
    \tilde{M}:=\frac12\sum_{l=1}^\infty\left|B_{2^{-(l-1)}}-B_{2^{-l}}\right|,
    \en
where $B$ denotes a Brownian bridge.
\bl\label{MNconvM}
There exists a constant $C>0$, such that
    \eq
    \Pbold\left(\left|M_N-\tilde M\right|>\eps\right)
    \le \frac{C N^2}{\eps}2^{-N/2}
    \en
for every $\eps>0$.
In particular, $M_N$ converges to $\tilde{M}$ in distribution as $N\to\infty$.
Moreover, it converges in $L^1$ and in $L^2$.
\el
We show in Proposition \ref{thm-1} below that $\tilde{M}\eqd M$.
The proof of Lemma \ref{MNconvM} is deferred to the appendix.
The combination of Lemmas \ref{MNformula} and \ref{MNconvM} yields the convergence
    \eq
    \frac1{\su\sqrt{n}}\;\max_{k=1,\dots,n}\left|I^{(B)}_{\text{out}_{k}}-k\,I_\text{lsb}\right|
    \longrightarrow\frac12\sum_{l=0}^\infty\left|B_{2^{-l}}-B_{2^{-(l+1)}}\right|
    \en
in distribution, in $L^1$ and $L^2$, as $N\to\infty$ (and thus also $n=2^N-1\to\infty$).
This proves \refeq{limitM}.


\subsection{A representation of $M$ in terms of i.i.d.\ standard normals}
In this section, we will prove the following representation formula, which expresses $M$ in terms of independent standard normal random variables.
\bp[Rewrite of $\tilde M$ in terms of standard normals]\label{thm-1}
    Let $Z_1,Z_2,\dots$ be a sequence of i.i.d.\ standard normal distributed random variables.
    Then, $\tilde{M}$ can be expressed as
    \eq\lbeq{reprformM}
    \tilde{M}\eqd\frac12\sum_{l=1}^\infty2^{-(l+1)/2}\left|\sum_{j=1}^{l-1}2^{-(l-j)}Z_j-Z_l\right|.
    \en
    In other words, $\tilde{M}$ in \refeq{defM} has the same distribution as $M$ in \refeq{formulaM}.
\ep

In order to obtain the limit law of $\tilde{M}$, we have made use of the representation  (B1)
of the Brownian bridge. Now we will primarily use (B2).
We will essentially use the following well-known property of Brownian motion:

\begin{lemma}[The conditional law of the middle point]
\label{lem-midpoint}
Let $\{W_s\}_{s=0}^{\infty}$ be a standard Brownian motion. Then, the distribution of $W_{\frac t2}$ conditional on $W_t=z$ is a normal distribution with mean $\frac z2$ and variance $\frac t{4}$.
\end{lemma}

\begin{lemma}[The distribution of $\{B_{2^{-l}}\}_{l=1}^{\infty}$]
\label{lem-B2l}
The distribution of $\{B_{2^{-l}}\}_{l=1}^{\infty}$ is equal to
    \eq
    B_{2^{-l}}=\sum_{j=1}^{l} 2^{-(j+1)/2-(l-j)} Z_j,
    \en
where $\{Z_j\}_{j=1}^\infty$ are independent and identically distributed standard normal random variables.
\end{lemma}

\proof The proof is by induction in $l$. For $l=1$, we use Lemma \ref{lem-midpoint} together with
the fact that $\{B_s\}_{s=0}^1$ is a Brownian motion conditioned on $B_1=0$. This implies
that the distribution of $B_{\frac 12}$ is equal to a normal random variable with mean 0
and variance $\frac 14$. Therefore,
    \eq
    B_{\frac 12}=\frac 12 Z_1,
    \en
where $Z_1$ is a standard normal distribution. This initializes the induction hypothesis.
To advance it, assume that
    \eq
    B_{2^{-(l-1)}}=\sum_{j=1}^{l-1} 2^{-(j+1)/2-(l-1-j)} Z_j.
    \en
Then, again using Lemma \ref{lem-midpoint}, the distribution of
$B_{2^{-l}}$ conditionally on $B_{2^{-(l-1)}}$ is a normal distribution with mean
$1/2\, B_{2^{-(l-1)}}$ and variance $2^{-(l+1)}$.
Therefore, if we denote $Z_l=2^{(l+1)/2}(B_{2^{-l}}-\frac 12 B_{2^{-(l-1)}})$, then $Z_l$
is a standard normal random variable independent of $B_{2^{-(l-1)}}$. As a consequence,
we have that
    \eq
    B_{2^{-l}}
    =2^{-(l+1)/2}Z_l +\frac 12 B_{2^{-(l-1)}}
    =\sum_{j=1}^{l} 2^{-(j+1)/2-(l-j)} Z_j,
    \en
where in the last step, we have used the induction hypothesis.
\qed
\vskip0.5cm

\noindent
\proof[Proof of Proposition \ref{thm-1}]
As a consequence of Lemma \ref{lem-B2l}, we obtain the identity
    \eq\lbeq{BI}
    B_{2^{-(l-1)}}-B_{2^{-l}}
    =\sum_{j=1}^{l-1}2^{-(j+1)/2-(l-j)}Z_j-2^{-(l+1)/2}Z_l
    \en
Thus, by \refeq{defM}, we obtain \refeq{reprformM}.
\qed

\subsection{The moments of $M$}
In this section, we identify the first two moments of $M$:
\bl[Moments of $M$]\label{MomentsOfM}
The expectation of $M$ is given by \refeq{expecM}, that is
    \eq
    {\mathbb E}[M]
    =\frac1{\sqrt{2\pi}} \sum_{l=1}^{\infty}\big(2^{-l}-2^{-2l}\big)^{1/2}
    \approx0.839792,
    \en
    and the variance by
    \eq\lbeq{formulaVar}
    \Var(M)
    =\frac1\pi\sum_{1\le l<k} (2^{-l}-2^{-2l})^{1/2}\,(2^{-k}-2^{-2k})^{1/2}\,
    \Big(\sqrt{1-\rho_{lk}^2}-1+\rho_{lk}\arcsin(\rho_{lk})\Big)
    +\frac{1}{6}\left(1-\frac2\pi\right)
    \en
with
    \eq
    \rho_{lk}=\frac{-2^{-(l+k)}}{\sqrt{2^{-l}-2^{-2l}}\,\sqrt{2^{-k}-2^{-2k}}}.
    \en
The variance of $M$ can be approximated by $\Var(M)\approx0.080066$.
\el
\vskip0.5cm

In the proof of Lemma \ref{MomentsOfM}, we will make use of the following
property of the bivariate normal distribution:

\begin{lemma}[The expected product of absolute values of normals]\label{lem-ExpectedProdAbsValueOfNormals}
Let $Y_1$, $Y_2$ be two standard normal random variables having a bivariate normal joint distribution with correlation coefficient $\rho$. Then
    \eq
    \lbeq{expprodav}
    {\mathbb E}[|Y_1||Y_2|]
    =\frac{2}{\pi}\sqrt{1-\rho^2}+
    \frac{2\rho}{\pi} \arctan\Big(\frac{\rho}{\sqrt{1-\rho^2}}\Big)
    =\frac{2}{\pi}\left(\sqrt{1-\rho^2}+\rho\,\arcsin(\rho)\right).
    \en
\end{lemma}
For a proof of Lemma \ref{lem-ExpectedProdAbsValueOfNormals} see e.g.\ \cite[Excercise 15.6]{Kendall77}.

\proof[Proof of Lemma \ref{MomentsOfM}.]
The expression for the mean is easily derived.
We note that
    \eq\lbeq{defNl}
    N_l:=B_{2^{-(l-1)}}-B_{2^{-l}}
    =\sum_{j=1}^{l-1}2^{-(j+1)/2-(l-j)}Z_j-2^{-(l+1)/2}Z_l
    \qquad l=1,2,\dots
    \en
is normally distributed with mean 0 and variance
    \eq\lbeq{defVl}
    v_l=\sum_{j=1}^{l}2^{-(j+1)-2(l-j)}=2^{-l}-2^{-2l},
    \en
and rewrite $M$ in \refeq{formulaM} as $M=\frac12\sum_{l=1}^\infty|N_l|$.
Since, for a normal random variable $Z$ with mean 0 and variance $\sigma^2$, we have that
    \eq
    {\mathbb E}[|Z|]=\sqrt{\frac{2\sigma^2}{\pi}},
    \en
the representation formula \refeq{reprformM} allows us to identify the mean of the random variable $M$ as
    \eq\lbeq{app4}
    {\mathbb E}[M]
    =\frac12\sum_{l=1}^\infty {\mathbb E}\big[|N_l|\big]
    =\frac1{\sqrt{2\pi}} \sum_{l=1}^{\infty}\sqrt{v_l}
    =\frac1{\sqrt{2\pi}} \sum_{l=1}^{\infty}\big(2^{-l}-2^{-2l}\big)^{1/2}
    \approx0.839792.
    \en

For the variance of $M$ we expand
    \eq\lbeq{varM}
    \Var(M)=\frac14\left(2\sum_{1\le l<k} \text{Cov}\big(|N_l|,|N_{k}|\big)
    +\sum_{l=1}^{\infty} \Var\big(|N_l|\big)\right).
    \en
The variance term is not too hard, as
    \eq
    \Var\big(|N_l|\big)
    =\big({\mathbb E}[N_l^2]-{\mathbb E}[|N_l|]^2\big)
    =\left(1-\frac2\pi\right)v_l,
    \en
and therefore,
    \eq
    \sum_{l=1}^{\infty} \Var\big(|N_l|\big)
    =\left(1-\frac2\pi\right)\sum_{l=1}^{\infty} (2^{-l}-2^{-2l})
    =\frac23\left(1-\frac2\pi\right).
    \en
Fix now $1\le l<k$.
Then
    \eqalign\nonumber
    \text{Cov}\big(N_l,N_{k}\big)
    &=\text{Cov}\left(\sum_{j=1}^{l-1}2^{-(j+1)/2-(l-j)}Z_j-2^{-(l+1)/2}Z_l,
                    \sum_{j=1}^{l}2^{-(j+1)/2-(l-j)}Z_j\right)\\
    &=\sum_{j=1}^{l-1}2^{-(j+1)-(l-j)-(k-j)}-2^{-(l+1)-(k-l)}
    =-2^{-(l+k)}.
    \enalign
Therefore, $\big(N_l,N_{k}\big)$ is a bivariate normal distribution
with mean $(0,0)$, variances $\big(v_l,v_k)$,
and correlation coefficient
    \eq\lbeq{defRho}
    \rho_{lk}=
    \frac{\text{Cov}(N_l,N_k)}{\sqrt{\text{Var}(N_l)}\sqrt{\text{Var}(N_k)}}
    =\frac{-2^{-(l+k)}}{\sqrt{2^{-l}-2^{-2l}}\,\sqrt{2^{-k}-2^{-2k}}}.
    \en
Denoting $(Y_l,Y_k)$ a bivariate normal distribution with means 0,
variances 1 and correlation coefficient $\rho_{lk}$,
we have
    \eqalign
    \text{Cov}\big(|N_l|,|N_{k}|\big)
    &=\sqrt{v_l\,v_k}\,\text{Cov}(|Y_l|,|Y_k|)
    \lbeq{CovNlNk}\\
    &=\sqrt{v_l\,v_k}\,\big({\mathbb E}[|Y_l| |Y_k|]-{\mathbb E}[|Y_l|]\,{\mathbb E}[|Y_k|]\big)\nn\\
    &=\sqrt{v_l\,v_k}\,\Big({\mathbb E}[|Y_l| |Y_k|]-\frac{2}{\pi}\Big).\nn
    \enalign
By Lemma \ref{lem-ExpectedProdAbsValueOfNormals}, as well as
\refeq{varM}--\refeq{CovNlNk}, we obtain
    \eq
    \lbeq{valueVar}
    \Var(M)=\frac1\pi\sum_{1\le l<k} \sqrt{v_l\,v_k}\, \Big(\sqrt{1-\rho_{lk}^2}-1+\rho_{lk}\arcsin(\rho_{lk})\Big)
    +\frac{1}{6}\left(1-\frac2\pi\right).
    \en
Using \refeq{defVl} and \refeq{defRho} we can approximate $\Var(M)$ numerically, which yields \refeq{varianceOfM}.
\qed

Having completed the proofs of \refeq{limitM}--\refeq{varianceOfM}, the proof of Theorem \ref{thm-ThermoResults} is complete.

\subsection{Approximating the density of $M$}
\label{sec-densityM}
We now derive a formula for the density of $M$. Without loss of generality, we may assume that $\bar I_u=0$ and $\sigma_u=1$, i.e., the $I_u$ are standard normal distributed.

By denoting $\Nfat=(N_1,N_2,\dots)$ and $\Zfat=(Z_1,Z_2,\dots)$, we rewrite \refeq{defNl} with the help of the (infinite) matrix $L$ as $\Nfat=-L\cdot\Zfat$, where
    \eq
    L_{jl}=
        \begin{cases}
        -2^{-(j+1)/2-(l-j)}\qquad&\text{if $j<l$};\\
        2^{-(l+1)/2}\qquad&\text{if $j=l$};\\
        0\qquad&\text{if $j>l$}.
        \end{cases}
    \en
Note that $L$ is a lower triangular matrix.

We will approximate the density of the infinite sum $M=1/2\sum_{l=1}^\infty|N_l|$ by the finite sum
    \eq
    M^{\sss (m)}=\frac12\sum_{l=1}^m|N_l|,\qquad m\in\Nbold.
    \en
Writing $\Nfat^{\sss (m)}=(N_1,\dots,N_m)$, $\Zfat^{\sss (m)}=(Z_1,\dots,Z_m)$, and $L^{\sss (m)}$ for the upper left $m\times m$ corner of the infinite matrix $L$, we have that $\Nfat^{\sss (m)}=-L^{\sss (m)}\cdot\Zfat^{\sss (m)}$. In particular, $\Nfat^{\sss (m)}$ is normally distributed with mean $(0,\dots,0)$ and covariance matrix $\Sigma^{\sss (m)}=L^{\sss (m)}\cdot (L^{\sss (m)})^T$, where we write $(L^{\sss (m)})^T$ for the transpose of the matrix $L^{\sss (m)}$. Note that
    \eq
    \big(\Sigma^{\sss (m)}\big)_{jl}=\begin{cases}
    -2^{-(j+l)}\qquad&\text{if $j\neq l$}\\
    2^{-l}-2^{-2l}\qquad&\text{if $j=l$.}
    \end{cases}
    \en

Given the mean and covariance matrix of a multivariate normal distribution, its density is known to be
\eq\lbeq{jointDensity}
    f_{\Nfat^{\sss (m)}}(\n)=\frac1{(2\pi)^{m/2}}\,\frac1{(\det\Sigma^{\sss (m)})^{1/2}}\;\exp\left\{-\frac12{\n\left(\Sigma^{\sss (m)}\right)^{-1}\n^T}\right\},
\en
where $\n=(n_1,\dots,n_m)\in\Rbold^m$.
We write $|\Nfat^{\sss (m)}|$ for the pointwise absolute value $(|N_1|,\dots,|N_m|)$ of the $m$-dimensional vector $\Nfat^{\sss (m)}$. Its density is given by
\eq\lbeq{jointDensityAbsoluteValue}
    f_{|\Nfat^{\sss (m)}|}(\n)=
\sum_{\sigma\in\{-1,1\}^m}
\frac1{(2\pi)^{m/2}}\,\frac1{(\det\Sigma^{\sss (m)})^{1/2}}\;\exp\left\{-\frac12{\left(\sigma \n)(\Sigma^{\sss (m)}\right)^{-1}(\sigma \n)^T}\right\},\qquad \n\in[0,\infty)^m,
\en
where $(\sigma \n)=(\sigma_1 n_1,\dots,\sigma_m n_m)$. See e.g., \cite[(6.3.20)]{Bain91}.

The determinant $\det\Sigma^{\sss (m)}$ and the inverse $(\Sigma^{\sss (m)})^{-1}$ of the covariance matrix are easy to compute since $L^{\sss(m)}$ is a triangular matrix. The results are stated in the following two lemmas.
\bl\label{detSigmaM}
For all $m\in\Nbold$, the determinant of $\Sigma^{\sss (m)}$ is
    \eq
    \det\Sigma^{\sss (m)}=2^{-m(m+3)/2}.
    \en
\el
\proof
First we note that
\eq
    \det\Sigma^{\sss (m)}
    =\big(\det L^{\sss (m)}\big)^2.
\en
Since $L^{\sss (m)}$ is a triangular matrix, its determinant is obtained by multiplying the entries on the diagonal, and hence
\eq
    \det\Sigma^{\sss (m)}
    =\Big(\prod_{j=1}^m-L_{jj}^{\sss(m)}\Big)^2
    =\Big(\prod_{j=1}^m2^{-(l+1)/2}\Big)^2
    =2^{-m(m+3)/2}.
\en
\qed\vspace{0.5cm}

\bl\label{lem-SigmaInv}
For all $m\in\Nbold$, the inverse of $\Sigma^{\sss (m)}$ is given by
    \eq
    \left(\Sigma^{\sss (m)}\right)^{-1}_{jl}=
    \begin{cases}
    (j\wedge l)\,2^{(j+l)/2-1}\qquad&\text{for $j\neq l$};\\
    (j+1)\,2^{j-1}\qquad&\text{for $j=l$}.
    \end{cases}
    \en
\el
\proof
Since $L^{\sss (m)}$ is a triangular matrix, it is easy to see that the inverse $(L^{\sss (m)})^{-1}$ is given by
    \eq
    (L^{\sss (m)})^{-1}_{jl}=
    \begin{cases}
        2^{(j-1)/2}\qquad&\text{if $j>l$};\\
        2^{(j+1)/2}\qquad&\text{if $j=l$};\\
        0\qquad&\text{if $j<l$}.
    \end{cases}
    \en
For $j<l$ we have that
    \eq
    \left(\Sigma^{\sss (m)}\right)^{-1}_{jl}
    =\sum_{k=1}^{j\wedge l}((L^{\sss (m)})^{-1})_{jk}\,((L^{\sss (m)})^{-1})_{lk}
    =\sum_{k=1}^j2^{(j-1)/2}2^{(l-1)/2}
    =j\,2^{(j+l)/2-1},
    \en
whereas for $j=l$,
    \eq
    \left(\Sigma^{\sss (m)}\right)^{-1}_{jj}
    =\sum_{k=1}^{j}((L^{\sss (m)})^{-1})_{jk}^2
    =\sum_{k=1}^{j-1}2^{(j-1)/2}+2^{(j+1)/2}
    =(j+1)\,2^{j-1}.
    \en
\qed\vspace{0.5cm}

In order to calculate the density of $M^{\sss (m)}=\frac12\big(|N_1|+ \dots+ |N_m| \big)$ at $y\ge0$, we have to integrate \refeq{jointDensity} over the $(m-1)$-dimensional surface
$\{(n_1,\dots,n_m)\,|\,2y=|n_1|+\dots+|n_m|\}$.
This leads to the following formula:
    \eqsplit\lbeq{densityFormula}
    f_{M^{\sss (m)}}(y)
    =&2\int\limits_0^{2y}\,\int\limits_0^{2y-n_1}\cdots\!\!\!\!\!\int\limits_0^{2y-n_1-\cdots-n_{m-2}}\!\!\!\!\!\!\!\!\!\!\!\!\!\!
    \quad f_{|\Nfat^{\sss (m)}|}\!\left(n_1,n_2,\dots,n_{m-1},(2y-n_1-\dots-n_{m-1})\right)\\
    &\hspace{5cm}
    dn_{m-1}\,\cdots\,dn_2\,dn_1,
    \qquad y\in[0,\infty).
    \ensplit
Note that there are $m-1$ integrals.
Thus for $m=1$ there are no integrals and
    \eq
    f_{M^{\sss (1)}}(y)=2f_{|\Nfat^{\sss (1)}|}(2y)=4\sqrt{\frac2\pi}\exp\{-8y^2\},
    \qquad y\in[0,\infty).
    \en

Finally, \refeq{densityFormula} gives us a formula for the density of
$M^{\sss (m)}$. The only (numerical) problem are the integrals.
\begin{figure}[h]
  \centering
  \includegraphics[width=13cm]{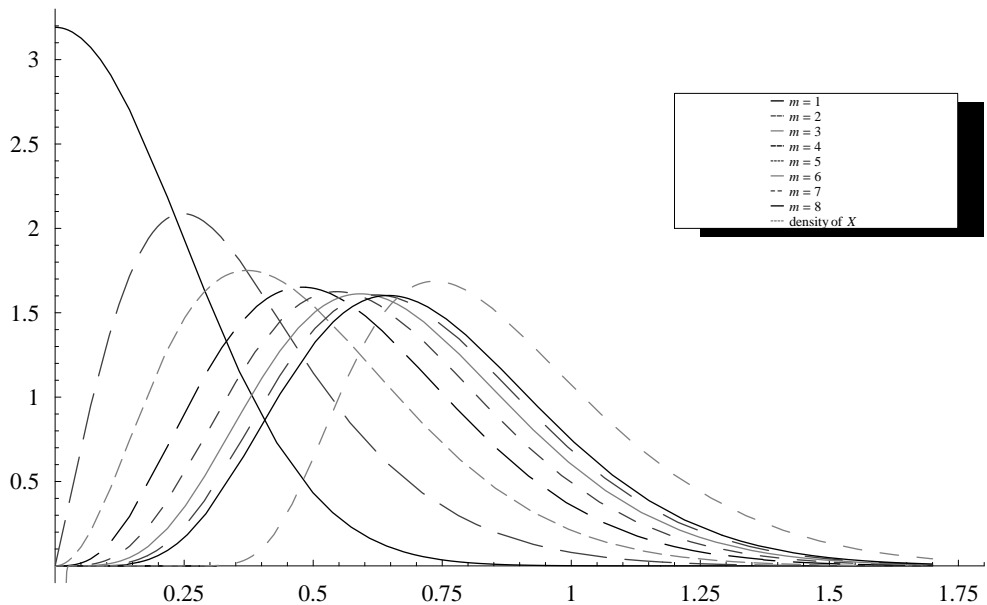}\\
  \caption{The density of $M^{\sss (m)}$ for different values of $m$, in contrast to the density of $X$.}
  \label{figureM0}
\end{figure}
For larger values of $m$, an intelligent way of numerical integration
seems necessary. However, for sufficiently small $m$, a mathematical
standard package, such as Mathematica, gives good approximations
(see Figure \ref{figureM0}).

\subsection{A disintegration approach to the density of $M$}\label{sectionDensityM2}
In this subsection we present a different approach to the density of $M$.
We define the quantity
\begin{equation}\lbeq{eqDefBarM}
    \M=\frac12\sum_{l=1}^\infty\left|W_{2^{-(l-1)}}-W_{2^{-l}}\right|
\end{equation}
for a Wiener process $\{W_s\}_{s\in[0,1]}$.
It is obvious that $M\eqd \M\big|_{W_1=0}$, cf.\ \refeq{defM} and (B2) on page \pageref{B2}.
Let $\hat f$ be the Fourier transform of the joint distribution of $\M$ and $W_1$, i.e.\
\begin{equation}
    \hat f(k_1,k_2)=\E\big[\exp\{i(k_1\M+k_2 W_1)\}\big].
\end{equation}
Using the independence and stationarity of the increments of the Wiener process in \refeq{eqDefBarM} we obtain
\begin{equation}
    \hat f(k_1,k_2)=\prod_{l=1}^\infty\underbrace{\E\big[\exp\{i(k_1\,|W_{2^{-(l-1)}}-W_{2^{-l}}|+k_2\,(W_{2^{-(l-1)}}-W_{2^{-l}}) )\}\big]}_{=\hat h_l(k_1,k_2)},
\end{equation}
where
\begin{equation}\label{}
    \hat h_l(k_1,k_2)=\E\big[\exp\{i(k_1|Z|+k_2 Z)\}\big],
\end{equation}
and $Z$ is a $\mathcal{N}(0,2^{-l})$-distributed normal random variable.
Once we have computed $\hat h_l$, we obtain the joint density of $(\M,W_1)$ via inverse Fourier transformation as
\begin{equation}
    f_{\M,W_1}(x,y)=\int\!\!\int\exp\{-ik_1x-ik_2y\}\,\hat f(k_1,k_2)\,\frac{dk_1\,dk_2}{(2\pi^2)^2}.
\end{equation}
The density of $M$ equals the density of $\M$ conditioned on $W_1=0$, whence
\begin{equation}\lbeq{newDensityM}
    f_M(x)=\frac{f_{\M,W_1}(x,0)}{f_{W_1}(0)}
    =\sqrt{2\pi}\int\!\!\int\exp\{-ik_1x\}\,\hat f(k_1,k_2)\,\frac{dk_1\,dk_2}{(2\pi^2)^2}.
\end{equation}

It remains to calculate $\hat h_l(k_1,k_2)$.
If we let
\begin{equation}\lbeq{defHatH}
    \hat h_l(k)=\E\big[\e^{ikZ}\1_{\{Z\ge0\}}\big],
\end{equation}
then $\hat h_l(k_1,k_2)=\hat h_l(k_1+k_2)+\hat h_l(k_1-k_2)$.
The dependence of $\hat h_l$ on $l$ can easily be eliminated by scaling:
\begin{equation}\label{}
    \hat h_l(k)=\hat h(2^{-l/2}k),
\end{equation}
where $\hat h$ is as in \refeq{defHatH} with $Z$ being $\mathcal{N}(0,1)$-distributed.
It can be shown that
\begin{equation}\label{}
    \hat h(k)=\frac12 \e^{-k^2/2} + i \frac{\e^{-k^2/2}}{2\sqrt{2\pi}}\int_{-k}^k\e^{x^2/2}\,dx,
\end{equation}
but the inverse Fourier transform in \refeq{newDensityM} seems intractable.

\section{Conclusions}
We have derived the distribution of the $\INLm$ in terms of Brownian
bridges. This distributional identity holds for all architectures, in particular also for
the segmented one.

In the thermometer case we have identified the limiting distribution
of $\INLm$ as the absolute maximum of Brownian bridges, which is
well-known in the literature.

For the binary case, we have identified the limiting distribution
of $\INLm$ in terms of a Brownian bridge. We have further
provided a representation in terms of
independent standard normal variables and have computed the mean and
variance of $\INLm$. Finally, we have given a procedure that approximates the
density.

We want to emphasize that the $\INLm$ in the thermometer case and
in the binary case behave differently. Although the densities
look alike, e.g.\ the upper tails are quite close to each other,
there are significant changes in the lower tail.
The thermometer
case has, compared to the binary case,
a slightly larger mean, but a slightly smaller variance.
Even though the distributions in the two cases are close, they are \emph{not} the same.

We still miss the distribution function for the binary case.
Random sums of the type in \refeq{defM} have appeared in the literature.
In particular the quantity
    \eq
    S=\sum_{l=1}^\infty 2^{-l} V_l
    \en
where $\{V_l\}_{l=1}^{\infty}$ are i.i.d.\ exponential random variables
arise in a variety of applications, see e.g.\ Ott, Kemperman and Mathis \cite{OKM96}, Guillemin, Robert and Zwart \cite{GZR04} and Litvak and van Zwet \cite{Litvak04}.
The density of $S$ can be expressed in terms of an infinite sum, cf.\ \cite[Section 5]{OKM96}.
When the summands are independent {\it uniform} random variables,
i.e., $\{V_l\}_{l=1}^{\infty}$ are i.i.d.\ uniform random variables on $(0,1)$,
the density of $S$ can be computed explicitly \cite{Fey06}.

Furthermore, it would be interesting to extend the results to the segmented case.
In particular, it would be of interest to investigate which limiting
$\INLm$ distribution has the smallest mean. This should correspond to the
optimal DAC architecture.
Practical implications of our results can be found in a
companion paper \cite{RaduHeydHofsHegtRoer06}.

\begin{appendix}
\section{Proof of Lemmas \ref{ConvergenceMaxB} and \ref{MNconvM}}
\proof[Proof of Lemma \ref{ConvergenceMaxB}.]
We first prove \refeq{thermo7}.
Let $\{B_s\}_{s\in[0,1]}$ be a Brownian bridge. Then,
\eq\lbeq{thermo2}
    \Big|\max_{k=1,\dots,n}\left|B_{k/n}\right|-\max_{t\in[0,1]}\left|B_{t}\right|\Big|
    \le
    \max_{\substack{k=1,\dots,n \\t\in[\frac{k-1}{n},\frac kn]}}
    \left|B_{k/n}-B_{t}\right|
\en
Using representation (B1) above, we can further bound \refeq{thermo2} from above by
\eq\lbeq{thermo3}
    \max_{\substack{k=1,\dots,n \\t\in[\frac{k-1}{n},\frac kn]}}
    \left|W_{k/n}-W_{t}\right|
    +\frac1n|W_1|
\en
for a Wiener process $\{W_s\}_{s\in[0,1]}$.
Using the Markov property and Brownian scaling, we obtain that for $k=1,\dots,n$,
\eq
    \max_{{t\in[\frac{k-1}{n},\frac kn]}}\left|W_{k/n}-W_{t}\right|
    \eqd \frac1{\sqrt{n}}\max_{t\in[0,1]}|W_t|,
\en
where $\eqd$ stands for equality in distribution.
Hence, for $C>0$,
\eqalign
    &\Pbold\!\left(\Big|\max_{k=1,\dots,n}\left|B_{k/n}\right|-\max_{t\in[0,1]}\left|B_{t}\right|\Big|
    \ge C\sqrt{\frac{\log n}{n}}\right)\nn\\
    &\le\Pbold\!\left(\max_{k=1,\dots,n}\bigg\{\max_{t\in[\frac{k-1}n,\frac kn]}\left|W_{k/n}-W_{t}\right|
    \ge \frac C2\sqrt{\frac{\log n}{n}}\bigg\}\right)
    +\Pbold\!\left(\frac1n|W_1|\ge\frac C2\sqrt{\frac{\log n}{n}}\right)\nn\\
    &\le n\,\Pbold\!\left(\max_{t\in[0,1]}\left|W_{t}\right|
    \ge \frac C2\sqrt{\log n}\right)
    +\Pbold\!\left(|W_1|\ge\frac C2\sqrt{n\log n}\right)\lbeq{thermo4}.
\enalign
For the first term, we bound for every $b\ge0$,
\eq
    \Pbold\!\left(\max_{t\in[0,1]}\left|W_{t}\right|\ge b\right)
    \le 2\,\Pbold\!\left(\max_{t\in[0,1]}W_{t}\ge b\right)
    = 4\,\Pbold\!\left(W_{1}\ge b\right)
    \le 4\e^{-b^2/2},
\en
where we use the reflection principle \cite[Theorem (6), p.\ 526]{Grimmett92}
in the second and a standard bound on the tail of standard normals in the third step.
Substituting $b=C/2\,\sqrt{\log n}$ we obtain
\eq\lbeq{thermo5}
    n\,\Pbold\!\left(\max_{t\in[0,1]}\left|W_{t}\right|
    \ge \frac C2\sqrt{\log n}\right)
    \le n\,4n^{-(C/2)^2/2}
    =4n^{1-C^2/8}.
\en
For the second term in \refeq{thermo4}, we obtain analogously
\eq\lbeq{thermo6}
    \Pbold\!\left(|W_1|\ge\frac C2\sqrt{n\log n}\right)
    \le 2\,\Pbold\!\left(W_1\ge\frac C2\sqrt{n\log n}\right)
    \le 2 n^{-C^2n/8}
\en
The bound \refeq{thermo4}, together with \refeq{thermo5}-\refeq{thermo6} proves \refeq{thermo7}.

For the convergence in $L^1$, we use that, with
    \eq
    X_n=\max_{k=1,\dots,n}\left|B_{k/n}\right|, \qquad X=\max_{t\in[0,1]}\left|B_{t}\right|,
    \en
we have that
    \eq
    \Ebold|X_n-X|=\int_{0}^\infty \Pbold(|X_n-X| \geq t)\,dt.
    \en
We split between $t\leq 4\sqrt{\frac{\log n}{n}}$ and $t>4\sqrt{\frac{\log n}{n}}$.
For $t\leq 4\sqrt{\frac{\log n}{n}}$, we bound $\Pbold(|X_n-X| \geq t)\leq 1$,
while for $t>4\sqrt{\frac{\log n}{n}}$, we use \refeq{thermo5} and \refeq{thermo6} to bound
    \eq
    \Pbold(|X_n-X| \geq t) \leq 6 n e^{-{nt^2}/{2}}.
    \en
Substitution of these bounds yields
    \eq
    \Ebold|X_n-X|=\int_{0}^\infty \Pbold(|X_n-X| \geq t)\,dt\leq
    4\sqrt{\frac{\log n}{n}}+6n\int_{4\sqrt{\frac{\log n}{n}}}^\infty e^{-{nt^2}/{2}}\,dt
    =O\Big(\sqrt{\frac{\log n}{n}}\Big).
    \en
The convergence in $L^2$ follows similarly, now using
    \eq
    \Ebold[(X_n-X)^2]=\int_{0}^\infty 2t\, \Pbold(|X_n-X| \geq t)\,dt.
    \en
We leave the details to the reader.
\qed
\vspace{.5cm}

\proof[Proof of Lemma \ref{MNconvM}.]
We observe that, by (\ref{eqM}), \refeq{defMN} and \refeq{defM},
and replacing $m$ by $N-m$,
    \eq
    \left|M_N-\tilde M\right|
    \le \sum_{m=1}^{N-1}\left|B_{2^{-m}}-B_{(2^{N-m}\,-1)/n}\right|
    +\sum_{m=N}^\infty\left|B_{2^{-m}}-B_{2^{-m-1}}\right|.
    \en
Therefore, for an arbitrary chosen constant $\eps>0$, we have that
    \eq\lbeq{app2}
    \Pbold\!\left(\left|M_N-\tilde M\right|>\eps\right)
    \le \underbrace{\Pbold\left(\sum_{m=1}^{N-1}\left|B_{2^{-m}}-B_{(2^{N-m}\,-1)/n}\right|>\frac{\eps}{2}\right)}_{(I)}
    +\underbrace{\Pbold\left(\sum_{m=N}^\infty\left|B_{2^{-m}}-B_{2^{-m-1}}\right|>\frac \eps2\right)}_{(II)}.
    \en
Using representation (B1), we see that for a Brownian bridge $\{B_s\}_{s\in[0,1]}$, and any $0\le s<t\le1$, the following holds for every constant $C>0$:
    \eqalign\nn
    \Pbold\left(|B_t-B_s|>C\right)
    &\le \Pbold\left(|W_t-W_s-(t-s)W_1|>C\right)\\
    &\le \Pbold\left(|W_t-W_s|>C/2\right)+\Pbold\left(|(t-s)W_1|>C/2\right).
    \enalign
Using the Markov inequality and Gaussian scaling we obtain
    \eq\lbeq{app1}
    \Pbold\left(|B_t-B_s|>C\right)
    \le \frac2{C}\sqrt{t-s}\,\,\E[|W_1|]+\frac2{C}(t-s)\,\E[|W_1|]
    \le\frac4{C}\sqrt{\frac2\pi (t-s)}.
    \en
We use \refeq{app1} to bound $(I)$ in \refeq{app2} from above by
    \eq\lbeq{app6}
    \sum_{m=1}^{N-1}\Pbold\!\left(\left|B_{2^{-m}}-B_{(2^{N-m}\,-1)/n}\right|>\frac{\eps}{2(N-1)}\right)
    \le \sum_{m=1}^{N-1}\frac{8(N-1)}{\eps}\sqrt{\frac1{2^m}-\frac{2^{N-m}\,-1}{2^{N}-1}}
    \le \frac{8N^2}{\eps\,2^{N/2}},
    \en
which converges to 0 as $N\to\infty$.
The second term $(II)$ in \refeq{app2} can be bounded using the Markov inequality by
    \eq
    (II)\le\frac2{\eps}\, \E\left[\sum_{m=N}^\infty\left|B_{2^{-m}}-B_{2^{-m-1}}\right|\right].
    \en
Using \refeq{app4}, this expectation can be computed as
    \eq
    \E\left[\frac12\sum_{m=N}^\infty\left|B_{2^{-m}}-B_{2^{-m-1}}\right|\right]
    =\frac{1}{\sqrt{2\pi}}\sum_{m=N}^\infty\left(2^{-m}-2^{-2m}\right)^{1/2},
    \en
and
    \eq\lbeq{app5}
    (II)\le
    \frac{2\sqrt{2}}{\eps\sqrt{\pi}}\sum_{m=N}^\infty\left(2^{-m}-2^{-2m}\right)^{1/2}
    \le\frac{4}{\eps\sqrt{\pi}(\sqrt 2-1)}\,2^{-N/2}
    \en
which converges to 0 as $N\to\infty$.

The combination of \refeq{app2}, \refeq{app6} and \refeq{app5} shows that, for $C=8+4(\sqrt{2\pi}-\sqrt\pi)^{-1}$,
    \eq
    \Pbold\left(\left|M_N-\tilde M\right|>\eps\right)
    \le \frac{C N^2}{\eps}2^{-N/2}
    \en
for every $\eps>0$, i.e., $M_N$ converges to $\tilde M$ in distribution as $N\to\infty$.

The convergence in $L^1$ and in $L^2$ follows as in the proof of Lemma \ref{ConvergenceMaxB}.
\qed
\end{appendix}

\subsection*{Acknowledgement}
The work of MH and RvdH was supported by the Netherlands Organisation for Scientific Research (NWO), and the work of GR was supported by STW, project {ECS.6098}.
We thank Marko Boon for help with the simulations in Figure \ref{figureM0}.
We thank David Brydges for enlightening discussions on multivariate normal distributions,  and Olaf Wittich for pointing our attention to the disintegration approach in Section \ref{sectionDensityM2}.

\bibliographystyle{abbrv}

\end{document}